\documentclass[12pt]{article}

\usepackage{amsmath, amsxtra, amsfonts,amscd,amssymb}
\usepackage{graphicx,wrapfig}
\usepackage{epstopdf}
\usepackage{url}
\usepackage[algo2e,linesnumbered, vlined,ruled]{algorithm2e}
\usepackage{float}
\usepackage{multirow}
\usepackage[margin=1in]{geometry}
\usepackage{color}  
\usepackage{mathrsfs}
\usepackage{comment}

\usepackage{tikz}
\usetikzlibrary{automata,arrows,positioning,calc}

\newtheorem{theorem}{Theorem}[section]
\newtheorem{prop}[theorem]{Proposition}

\newtheorem{corollary}[theorem]{Corollary}

\def\d{\mathrm{d}}

\def\A{\mathcal{A}}

\title{\LARGE \bf
Bilinear Controllability of a Simple Reparable System
}

\author{Daniel Owusu Adu$^{1}$ and Weiwei Hu$^{2}$
\thanks{$^{1}$Daniel Owusu Adu is with the Department of Mathematics, University of Georgia, Athens, GA, USA.
        {\tt\small  daniel.adu@uga.edu}}%
\thanks{$^{2}$Weiwei Hu is with the Department of Mathematics, University of Georgia, Athens, GA, USA. %
        {\tt\small Weiwei.Hu@uga.edu. }The author acknowledges financial support from NSF grant no. DMS-2229345, 2205117, 2111486
 and AFOSR, FA9550-23-1-0675.}%
}

\begin{document}

\maketitle
\thispagestyle{empty}
\pagestyle{empty}

\begin{abstract}

 Reparable systems are systems that are characterized by their ability to undergo maintenance actions when failures occur. These systems are often described by transport equations, all coupled through an integro-differential equation. In this paper,  we address the understudied aspect of the controllability of reparable systems. In particular, we focus on a two-state reparable system and our goal is to design a control strategy that enhances the system availability- the probability of being operational when needed. We establish bilinear controllability, demonstrating that appropriate control actions can manipulate system dynamics to achieve desired availability levels. We provide theoretical foundations and develop control strategies that leverage the bilinear structure of the equations.

\end{abstract}

\section{Introduction}
In the realm of reliability and maintenance, reparable systems play a pivotal role, where maintenance strategies are carefully devised to promptly address failures and minimize total breakdown. Reparable systems occur naturally in problems of product design, inventory systems, computer networking, electrical power system  and complex manufacturing processes.   These systems are characterized by their ability to undergo repair actions when failures occur, ensuring their continuous functionality. Extensive research has been conducted in this domain, particularly focusing on systems with arbitrarily distributed repair times, often governed by complex systems of coupled partial and integro-differential  hybrid equations (e.g.,\,\cite{Chung-1,  gao2022stability,  Guo-1, gupur2011functional, HR-2, HXYZ-1,hu2007exponential, xu2005asymptotic}).

The primary emphasis in much of this prior research was concerned with the well-posedness of mathematical models and its the asymptotic behavior. While these endeavors have contributed significantly to our understanding of the reparable systems, our current work takes a distinct perspective. We seek to shed light on the critical issue of controllability concerning maintenance strategies with a specific aim: enhancing the availability of reparable systems with prescribed demands. 

Availability, in this context, refers to the probability that a reparable system remains operational and free from failure or undergoing repair actions when it is needed for use \cite{jardine2013maintenance, moubray1997reliability}. To address this concern, this work mainly focuses on a particular class of reparable multi-state systems with arbitrarily distributed repair time, initially introduced by Chung \cite{Chung-1}. The system is characterized by coupled transport and integro-differential equations, which proposing a challenging yet realistic model for various real-world applications.

In \cite{Chung-1}, it is assumed that there are $M$ distinct failure modes associated with a device and initially, the device is in the good mode, denoted by $0$. Transitions between states $0$ and $j$ are allowed, with $j=1,2, \dots, M$, which are determined by failure rates and repair rates. Repair times, on the other hand, follow arbitrary distributions, adding an element of unpredictability to the maintenance process.
The repair actions undertaken in our model are naturally likened to corrective maintenance, which serves to restore a failed system to operational
status. In fact, optimal repair rate design of such systems has been discussed in our previous work \cite{boardman2019optimal, hu2022optimal}. It gave rise to a bilinear open-loop control problem. In contrast,  our current is aimed at constructing a space-time dependent repair rate design for achieving exact controllability of the system under certain conditions.

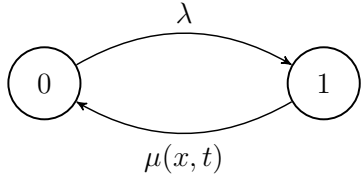
\begin{figure}
\centering
\resizebox{0.3\textwidth}{!}{
  \begin{tikzpicture}[->, >=stealth', auto, semithick, node distance=4cm]
    \tikzstyle{every state}=[fill=white, draw=black, thick, text=black]
    \node[state] (A) {$0$};
    \node[state] (B) [right of=A] {$1$};
    \path
    (A) edge[bend left, above] node{$\lambda$} (B)
    (B) edge[bend left, below] node{$\mu(x,t)$} (A);
  \end{tikzpicture}
}
\caption{Transition diagram of the reparable two-state system}
\label{transition}
\end{figure}

In this work, we focus on a simplified scenario where the reparable system has only one failure mode and  the repair rate is  allowed to depend on the system running time which is more realistic (see Fig.~\ref{transition}). While  one failure mode  may seem restrictive, it captures the essence of the original model and serves the purpose to convey the core principles of our maintenance strategy design.

The precise model characterizing the two-state reparable system reads
\begin{eqnarray}
&\displaystyle\frac{d p_0(t)}{d t} =-\lambda  p_0(t)+ \int_0^{L} \mu(x)p_{1}(x,t) \,d x,   \label{eqn:DistParaSys1}\\
&\displaystyle\frac{\partial{p_{1}(x,t)}}{\partial{t}} + \frac{\partial{p_{1}(x,t)}}{\partial{x}} =-\mu(x)p_{1}(x,t),
 \label{eqn:DistParaSys2}
\end{eqnarray}
with boundary condition
\begin{equation}
p_{1}(0,t)=\lambda p_0(t),  
\label{eqn:DistParaSys_BC}
\end{equation}
and  nonnegative initial conditions
\begin{equation}
p_0(0)=p_{0_0}, \quad p_{1}(x,0)=p_{1_0}(x).
\label{eqn:DistParaSys_IC}
\end{equation}
Here
\begin{enumerate}
 \item $p_0 (t)$:  probability that the device is in good mode $0$ at time $t>0$;
 \item $p_{1}(x, t)$: probability density (with respect to repair time $x$) that the failed device is in failure mode $1$ at time $t>0$ and has an elapsed repair time of
  $x\in [0, L]$ for $0<L<\infty$.
   Let  $\hat{p}_{1}(t)$  denote the probability of the  device in failure
mode $1$ at time $t$, then $\hat{p}_{1}(t)$ is given  by 
\begin{align}
            \hat{p}_1(t)=\int^{L}_{0}p_{1}(x,t)\,d x; \label{probability_p1}
\end{align}

  \item $\lambda>0$: constant failure rate of the device for failure mode $1$;
  \item $\mu(x)\geq0$: repair rate of the device  with an elapsed repair time  $x$.  Assume that for $0<l<L$,
\begin{align}
\int^{l}_{0} \mu(x)\,dx <\infty\quad\text{and}\quad \int^{L}_{0} \mu(x)\,dx =\infty.
\label{mu}
\end{align}
\item 
 The initial probability distributions of the system in good and failure modes  satisfy 
\begin{align}\label{EST_ini}
 p_{0_0}+\int^{L}_0p_{1_0}(x)\, dx=1. 
 \end{align}
   \end{enumerate}
The following assumptions are associated with the device:
\begin{enumerate}
\item The failure rates are constant;
\item  All failures are statistically independent;
\item  All repair time of failed devices are arbitrarily distributed;
\item There is only one mode of failure denoted by $1$, and $0$ implies the good state;
\item Repair is to like-new and it does not cause damage to any other part of the system.
\item Transitions are permitted only between states 0 and 1;
\item The repair process begins soon after the device is in failure state;
\item The repaired device is as good as new;
\item  No further failure can occur when the device has been down.
\end{enumerate}
\subsection{Well-posedness and stability of the Model }

The well-posedness and stability issues of system \eqref{eqn:DistParaSys1}--\eqref{eqn:DistParaSys_IC} for given failure and time-independent repair rates have been  well studied  in  \cite{HXYZ-1, hu2016differentiability, xu2005asymptotic} by using $C_0$-semigroup theory in a nonreflexive Banach space $X=\mathbb{R}\times L^1(0, L)$ equipped with the norm $\|\cdot\|_{X}=|\cdot|+\|\cdot\|_{L^1}$.  It is proven in \cite{xu2005asymptotic} that system operator  generates a positive $C_0$-semigroup of contraction. Thus the solution to \eqref{eqn:DistParaSys1}--\eqref{eqn:DistParaSys_IC} is nonnegative if the initial data are nonnegative.  In fact, using the method of characteristics we get
\begin{align*}
p_0(t)=p_{0_0}e^{-\lambda t} + \int^{t}_{0}e^{-\lambda(t-\tau) } \int^{\tau}_{0}\mu(x)p_{1}(x,\tau) dx d\tau.
\end{align*}
and 
\begin{align*}
&p_1(x,t)=\left\{\begin{array}{c}
\displaystyle 
\lambda p_{0}(t-x)e^{-\int^{x}_{0}\mu(s) ds}, \quad x<t, \\ 
\displaystyle
  p_{1_0}(x-t)e^{-\int^{x}_{x-t}\mu(s)\,ds}, \quad  x\geq t.\\
 \end{array}\right. 
\end{align*}
It is easy to verify that $p_0\in W^{1,\infty}(0, \infty)$ and $p_1\in L^\infty(0, \infty; W^{1,1}(0, L))$ for $t>x$, which is independent the regularity of the initial data.  For $x\geq t$, the same regularity results hold when $p_{1_0}\in W^{1,1}(0, L)$. Moreover,
\begin{align}
\int^{L}_{0} \mu(x)\,dx =\infty\implies p_1(L, t)=0, \label{prop_p1}
\end{align}
 which indicates that the probability density distribution of the system in failure mode becomes zero once the repair time reaches its maximum.  As a result, one can show that $\frac{dp_0}{dt}+\frac{\partial \int^L_0p_1(x,t)\,dx}{\partial t}=0$, and hence from~\eqref{EST_ini}
\begin{align}
p_0(t)+\int^{L}_0p_{1}(x, t)\, dx=1, \label{sys_conser}
 \end{align}
for $\forall t>0$.
Furthermore, it can be shown that zero is a simple eigenvalue of the system operator and also a unique spectrum on the imaginary axis. The $C_0$-semigroup generated by the system operator is eventually compact. As a result, the time-dependent solution 
exponentially converges to the its steady-state solution, which is the eigenfunction associated with the zero eigenvalue given by 
\begin{align}
P_{ss}(x)=\left(p_{0_{ss}},  \lambda p_{0_{ss}}e^{-\int^{x}_{0} \mu(s)\,ds}\right)^{T}.
\label{ss}
\end{align}
The detailed proof can be found in    \cite{HXYZ-1, hu2016differentiability}.

   The present work will investigate the bilinear controllability of the system via a space-time dependent repair rate $\mu(x, t)$.  Due to the properties of nonnegativity and conservation of the system, the desired states should also satisfy these attributes described by \eqref{prop_p1}--\eqref{sys_conser} and the boundary condition \eqref{eqn:DistParaSys_BC}. 
 Moreover, we assume that the desired probability density distribution of the system in failure mode is a strictly decreasing function of the same regularity of the system solution. In other words, while under repair, it is not expected that the desired density distribution of the failure rate increases. 
 
\subsection{Problem Statment}
Given $t_f>0$ and a nonnegative initial datum   $ P_0(x)=(p_{0_0}, p_{1_0}(x))^T\in X$ satisfying \eqref{EST_ini}, let   $P^*(x)=(p^*_0, p^*_1(x))^T\in X$ be a desired nonnegative distribution satisfying $\|P_0\|_{X}=\|P^*\|_{X}=1$,
that is,
\begin{align}
p_0(0)+\int^L_0p_1(x, 0)\d x=p^*_0+\int^L_0p^*_1(x)\d x=1. \label{2sys_conser}
\end{align}
Moreover,  assume 
\begin{align}
&p^*_1(0)=\lambda p^*_0,  \ p^*_1(L)=0,\ p^*_1\in W^{1,1}(0, L), \label{1prop_desire_p1}\\
&p^*_1(x)>0\quad  \text{and}\quad  \frac{d p^*_1(x)}{dx}<0, \ \forall x\in (0, L). \label{2prop_desire_p1}
\end{align}
 Determine whether there exists  a space-time dependent repair rate $\mu(x,t)$ such that the solution  $P(x, t)=(p_0(t), p_1(x,t))^T$ to \eqref{eqn:DistParaSys1}--\eqref{eqn:DistParaSys2}    satisfies $P(\cdot, t_f)=P^*(\cdot).$

\section{Bilinear Control Design}
First of all,  observe that if the repair rate is time-independent, we set the steady-state solution in  \eqref{ss}  to be the desired distribution 
$(p^*_0, p^*_1(x))$, i.e., $p_{0_{ss}}=p^*_0$ and 
$ \lambda p_{0_{ss}}e^{-\int^{x}_{0} \mu(s)\,ds} =p^*_1(x).$
We obtain 
\begin{align}
\mu(x)=(-\ln p^*_1(x))'= -\frac{p^*_{1_x}}{p^*_1},  \label{prop_mu}
\end{align}
which  satisfies \eqref{mu}.
This also implies that if the repair rate $\mu(x)$ satisfies  \eqref{prop_mu}, then  the system solution converges to $(p^*_0, p^*_1(x))$ exponentially.

Next  we investigate the bilinear controllability of the repair rate when it is allowed to depend on system running time $t$. Note that for any $t_f>0$ we can always choose a constant $c_0>0$ such that 
  $c_0\sum^{\infty}_{k=1} \frac{1}{k^2}=t_f$.  In fact,  $\sum^{\infty}_{k=1} \frac{1}{k^2}=\frac{\pi^2}{6}$ and hence $c_0=\frac{6t_f}{\pi^2}$. 
 Inspired by \cite{elamvazhuthi2017controllability, elamvazhuthi2018bilinear}, we consider
 the space-time dependent repair rate $\mu(x, t)$   given by 
\begin{align}
\mu(x,t)=-\frac{p_{1_x}}{p_1}+\alpha i\frac{(g(x)p_1)_x}{p_1},  \label{mu_feedback}
\end{align}
for $ t\in [c_0\sum^{i-1}_{k=1} \frac{1}{k^2}, c_0\sum^{i}_{k=1} \frac{1}{k^2})$ and $i\in \mathbb{Z}^{+}$, where $g(x)=\frac{1}{p^*_1(x)}$
  and $\alpha >0$  is a constant to be properly chosen.   Here   we set $ \sum^{i}_{k=1} \frac{1}{k^2}=0$ if $i=0$ and let $t_i=c_0\sum^{i}_{k=1} \frac{1}{k^2}, i\in \mathbb{Z}^+$, in the rest of our discussion.  In this case, $ \mu$ is piecewise defined in time $t$
  and it is straightforward to verify that 
  $\mu$ satisfies \eqref{mu} if $\alpha \geq \sup_{x\in [0, L]}|p^*_1(x)|=p^*_1(0)$. 

Our  main result  of this work is stated as follows. 

\begin{theorem}\label{main}
Given a nonnegative  initial datum     $ P_0=(p_{0_0}, p_{1_0})^T$ with 
$\|P_0\|_{X}=1$,
for  any nonnegative  $P^*=(p^*_0, p^*_1)^T\in X$   satisfying conditions 
\eqref{2sys_conser}--\eqref{2prop_desire_p1}, 
there exists a space-time dependent repair rate $\mu$ defined by \eqref{mu_feedback},
such that  the  solution 
$P(\cdot,t)=(p_0(t), p_1(\cdot,t))^T$ to \eqref{eqn:DistParaSys1}--\eqref{eqn:DistParaSys_IC} satisfies $P^*(\cdot, t_f)=(p^*_0, p^*_1(\cdot ))^T$.

\end{theorem}

The proof of  Theorem \ref{main} contains several components.  We first establish the well-posedness and stability analysis of the closed-loop system.

\subsection{Well-posedness of the Closed-Loop System }
Replacing $\mu$ by \eqref{mu_feedback} in \eqref{eqn:DistParaSys1}--\eqref{eqn:DistParaSys2}  leads to the following closed-loop system
\begin{align}
&\displaystyle\frac{d p_{0,i}(t)}{d t}=\alpha i \int^L_0\frac{\partial{(g(x)p_{1,i})}}{\partial x} \d x,
 \label{1closed_loop}\\
& \frac{\partial{p_{1,i}(x,t)}}{\partial{t}} 
=-\alpha i \frac{\partial{(g(x)p_{1,i})}}{\partial x},
\ \label{2closed_loop}
\end{align}
where  $p_{0,i}(t):=p_{0}(t)$ and $p_{1,i}(x,t):=p_1(x,t)$ for $(x,t)\in(0,L)\times(t_{i-1}, t_i)$ and $i\in \mathbb{Z}^+$,  
with  boundary condition 
\begin{align}
p_{1,i}(0,t)=\lambda p_{0,i}(t) 
\ \label{IC_F} 
\end{align}
and the initial conditions 
\begin{align}
 &p_{0,i}(0)=p_{0,i-1}(t_{i-1}),\label{IC_F1}  \\
&p_{1,i}(x,0)=p_{1,i-1}(x, t_{i-1}).
 \label{IC_F2} 
\end{align}

We first solve the solution to the closed-loop system \eqref{1closed_loop}--\eqref{2closed_loop} using the method of characteristics.

Let $\phi_{0, i}(t)=p_{0,i}(t)$ and $\phi_{1,i}(x,t)=\alpha i g(x)p_{1,i}(x, t)$ for $ t\in (t_{i-1}, t_i)$ and $i\in \mathbb{Z}^+$. Then  
 \eqref{1closed_loop}--\eqref{2closed_loop} become
\begin{align}
&\displaystyle\frac{d \phi_{0,i}(t)}{d t}= \int_0^{L} \frac{\partial{\phi_{1,i}}}{\partial x} \d x,
 \label{1closed_loop_phi}\\
&\frac{\partial{\phi_{1,i}}}{\partial{t}}  =-\alpha i g (x) \frac{\partial \phi_{1,i}}{\partial x}, \label{2closed_loop_phi}
\end{align}
with boundary condition
\begin{align}
\phi_{1,i}(0,t)=\alpha ig(0)p_{1,i}(0,t)= \alpha ig(0)\lambda \phi_{0,i}(t)\label{closed_loop_BC_phi}
\end{align}
and initial conditions
\begin{align}
\phi_{0,i}(0)=p_{0,i}(0), \quad
 \phi_{1,i}(x,0)=\alpha i g(x)p_{1,i}(x, 0).
\label{closed_loop_IC_phi}
\end{align}

 Let  $\frac{dx}{dt}=\alpha i g(x), \quad x(0)=x_0.$
Then
 $\frac{dx}{\alpha i g(x)}=\frac{1}{\alpha i }p^*_1(x)\d x = dt.$
Let $\tilde{p}^*_{1,i}(x)=\frac{1}{\alpha i }\int^x_0 p^*_1(s)\,ds$. Then  $\tilde{p}^*_{1,i}(x)= t+\frac{1}{\alpha i }\int^{x_0}_0 p^*_1(s)\,ds$. Since
$ \frac{ d\tilde{p}^*_{1,i}}{dx}=\frac{1}{\alpha i } p^*_1>0$ for $ x\in (0, L)$ 
 by \eqref{2prop_desire_p1}, this  implies that  $\tilde{p}^*_{1,i}(x)$ is  a strictly monotonically  increasing  function for $ x\in [0, L]$, and hence invertible. 
 Let $\xi=\tilde{p}^*_{1,i}(x)- t$. Then $x=(\tilde{p}^*_{1,i})^{-1} (\xi+ t)$. 
 Define $\Psi_{1,i}(t)=\phi_{1,i}((\tilde{p}^*_{1,1})^{-1}(  t+\xi), t)$.  Then 
\begin{align}
\frac{d\Psi_{1,i}}{dt}&= \alpha i g(x) \frac{\partial \phi_{1,i}}{\partial x}+\frac{\partial \phi_{1,i}}{\partial t}
=0.\label{1sol_y1}
\end{align}

For  $\xi<0$, i.e., $\tilde{p}_{1,i}(x)<  t$, the solution is determined by the boundary condition, so we integrate \eqref{1sol_y1} from some $t$ such that $x=(\tilde{p}^{*}_1)^{-1}(\xi+ t)=0$, i.e., $\xi+  t=0$, and hence $t=-\xi$. Integrating   \eqref{1sol_y1}  from 
$-\xi$ to $t$  follows
\begin{align}
\Psi_{1,i}(t)&=\phi_{1,i}(0, -\xi)=\phi_{1,i}(0, -\tilde{p}^*_{1,i}(x)+ t))
\nonumber\\
&=\alpha i g(0)\lambda\phi_{0,i}\left(  t-\tilde{p}^*_{1,i}(x)\right). \label{1EST_phi1}
\end{align}

For  $\xi\geq 0$, i.e., $\tilde{p}^*_{1,i}(x)\geq   t$, the solution is determined by the initial  condition. So we integrate  \eqref{1sol_y1}  from $0$ to $t$  and obtain 
\begin{align}
\Psi_{1,i}(t)&=\phi_{1,i}((\tilde{p}^*_{1,i})^{-1}(\xi), 0)\nonumber\\
&=\phi_{1,i}((\tilde{p}^*_{1,i})^{-1}(\tilde{p}^*_{1,i}(x)-  t), 0). \label{2EST_phi1}
\end{align}
To simplify the notation, we let $$\psi(x, t)=(\tilde{p}^*_{1,i})^{-1}(\tilde{p}^*_{1,i}(x)- t), \quad  t_{i-1}<t\leq \tilde{p}^*_{1,i}(x). $$
 Therefore,
\begin{align*}
\phi_{1,i}(x, t)=
\left\{\begin{array}{c}
\displaystyle 
\phi_{1,i}(\psi(x, t), 0), \quad t\leq \tilde{p}^*_{1,i}(x) ;\\
\alpha i g(0)\lambda\phi_{0,i}(  t-\tilde{p}^*_{1,i}(x)), \ t>\tilde{p}^*_{1,i}(x).
 \end{array}\right. 
\end{align*}
Solving $\phi_{0,i}$ from \eqref{1closed_loop_phi} yields 

\begin{align*}
&\phi_{0,i}(t)=\int^t_0\big(\phi_{1,i}(L,  \tau)-\phi_{1,i}(0,  \tau)\big)\,d\tau+\phi_{0,i}(0).
\end{align*}

Therefore, for any initial datum $P_0=(p_{0_0}, p_{1,0})^T\in X$, if the solution $P_i(x,t)=(p_{0,i}(t), p_{1,i}(x,t))^T$ to  \eqref{1closed_loop}--\eqref{IC_F2} exists,  then it is given by 
\begin{align}
p_{0,i}(t)
=\left\{\begin{array}{c}
\displaystyle 
\alpha i \int^t_0 \big(g(\psi(L, \tau))p_1(\psi(L, \tau), 0)\\
-g(0)\lambda p_0(\tau)\big)\,d\tau+p_{0,i}(0), 
\ t\leq \tilde{p}^*_{1,i}(L) ;
\\
\alpha ig(0)\lambda\int^t_ 0\big(p_0(  \tau- \tilde{p}^*_{1,i}(L))
-p_0(\tau)\big)\, d\tau\\
+p_{0,i}(0), \ t>\tilde{p}^*_{1,i}(L),
 \end{array}\right.  
 \label{ex_sol_p0}
\end{align}
and
\begin{align}
p_{1,i}(x, t)= 
\left\{\begin{array}{c}
\displaystyle 
\frac{1}{g(x)}g(\psi(x, t))p_1(\psi(x, t), 0),\\ \quad t\leq \tilde{p}^*_{1,i}(x) ;\\
\frac{1}{g(x)} g(0)\lambda p_0(  t-\tilde{p}^*_{1,i}(x)),\\ \quad t>\tilde{p}^*_{1,i}(x).
 \end{array}\right. 
  \label{ex_sol_p1}
\end{align}

To focus on our discussion, we first assume that $0<p^*_1(0)\leq 1$ and   investigate the properties of the solution when  $\alpha=1$ and $i=1$. In this case,
$\tilde{p}^*_{1,1}(x)=\int^x_0 p^*_1(s)\,ds$. We can rewrite  \eqref{1closed_loop}--\eqref{IC_F2} as an abstract Cauchy problem in  $X$
\begin{align}
\left\{\begin{array}{ll}
\dot{P}(t)=\mathcal{A}P(t), \quad t>0,\\
P(0)=(p_{0_0}, p_{1_0})^{T},
\end{array}
\right. \label{IVP_closed_loop}
\end{align}
where $\mathcal{A}$ is defined by
\begin{align}
\mathcal{A}= \displaystyle 
\begin{pmatrix}
0 &  \int^L_0\frac{\partial{(g(x)\cdot)}}{\partial x} \d x \\ 
0  & - \frac{\partial{(g(x)\cdot)}}{\partial x} 
\end{pmatrix}
 \label{A_new}
\end{align}with domain
\begin{align*}
D(\mathcal{A})=&\big\{ P\in X\big| gp_1\in W^{1,1}(0, L)
\\ &\quad\mathrm{and}\ p_1(0)=\lambda p_0 \big\}.
\end{align*}

Note that  $gp_1\in W^{1,1}(0, L)$ implies  $gp_1\in C[0, L]$ by Sobolev imbedding, 
and hence
\begin{equation}\label{eq: finite at the end point}
\int^L_0\frac{\partial (g(x)p_1(x))}{\partial x}\,dx=g(L)p_1(L)-g(0)p_1(0)<\infty.
\end{equation}
Moreover,  since $\lim_{x\to L}g(x)=\lim_{x\to L}\frac{1}{p^*_1(x)}=\infty$ based on \eqref{1prop_desire_p1}, we must have $p_1(L)=0$.

\begin{theorem}
The system operator $\mathcal{A}$ defined by \eqref{A_new} generates a positive  $C_{0}$-semigroup of contraction, denoted by  $T(t), t\geq 0$, on $X$.
\end{theorem}

The proof follows the similar  procedure as in  that of \cite[Thm.\,2.1]{xu2005asymptotic}. The details are omitted  here due to space limit. 
\subsection{Exponential Stability }
To establish the exponential  stability of \eqref{IVP_closed_loop}, we first  note that   when $\mu(x,t)=-\frac{p_{1_x}}{p_1}+\frac{(g(x)p_1)_x}{p_1}$,  the desired distribution 
$P^*(x)=(p^*_0, p^*_1(x))^T $\label{ss_closed_loop}
is the only steady-state solution of the closed-loop system \eqref{IVP_closed_loop} and it is the eigenfunction associated with eigenvalue zero of $\mathcal{A}$. One can further show that  zero is a simple eigenvalue  and the only spectrum on the imaginary axis following the proof of \cite[Prop.\,2.1]{hu2016differentiability}.  In addition, we can show that $T(t)$ is eventually compact. 
\begin{theorem}
 The $C_0$-semigroup $T(t)$ is compact when $t>\max\{L, \|p^*_1\|_{L^1}\}$.
\end{theorem}
\noindent{\it Proof.} First of all, we can show that the resolvent operator  $\mathcal{R}(r, \A)$ is compact  for any $r\in \rho(\A)$ as in   \cite{hu2016differentiability}. According to \cite[Cor. 3.4, p.\,50]{pazy1983semigoroups},  it suffices to show that $T(t)$ is continuous in the uniform operator topology for  $t>\max\{L, \|p^*_1\|_{L^1}\}$, that is,
\[
\sup_{\|P_0\|_{X}\leq 1, P_0\neq 0}\|T(t+h)P_0-T(t)P_0\|_{X}\xrightarrow{\text{$h\rightarrow 0$}}0, 
\]
uniformly,  
where
 \begin{align*}
&\|T(t+h)P_0-T(t)P_0\|_{X}=|p_0(t+h)-p_0(t)|\\
&\qquad+\int^L_0|p_1(x, t+h)-p_1(x, t)|\,dx.
\end{align*}
Since $ \tilde{p}^*_{1,1}(x)$ is strictly monotonically  increasing for $x\in [0, L]$,  we have $\|p^*_1\|_{L^1}=\tilde{p}^*_{1,1}(L)\geq \tilde{p}^*_{1,1}(x)$.
Thus  when $t >\max\{L, \|p^*_1\|_{L^1}\}$, we get $t>x$ and $t>\tilde{p}^*_{1,1}(x)$.
With the help of \eqref{ex_sol_p0}--\eqref{ex_sol_p1}, we obtain 
\begin{align}
&|p_0(t+h)-p_0(t)| \nonumber \\
&\qquad=\left|\alpha g(0)\lambda \int^{t+h}_ t\big(p_0(  \tau-\tilde{p}^*_{1,i}(L))-p_0(\tau)\big)\,d\tau\right|
\nonumber\\
&\qquad \leq 2\alpha g(0)\lambda h \sup_{t\geq 0 } |p_0(t)|\label{eq: estimate of the good state}
\end{align}
and 
\begin{align}
&\int^L_0|p_1(x, t+h)-p_1(x, t)|\,dx \nonumber\\
&=\int^L_0\Big|\frac{1}{g(x)} g(0)\lambda\big (p_0(  t+h-\tilde{p}^*_{1,1}(x))\nonumber\\
&\qquad-p_0(  t-\tilde{p}^*_{1,1}(x))\big)\Big|\,dx\nonumber\\
&=g(0)\lambda\int^L_{0}\Big|p_0(  t+h-\tilde{p}^*_{1,1}(x))\nonumber\\
&\qquad-p_0(  t-\tilde{p}^*_{1,1}(x))\Big|p_1^*(x)\,dx.
\label{Int_p0}
\end{align}
Let
$
\tilde{t}=t-\tilde{p}^*_{1,1}(x)>0,
$
then $
d\tilde{t}=-p^{*}_1(x)dx,
$
and hence \eqref{Int_p0} becomes 
\begin{align}
&g(0)\lambda\int^{L}_0\Big|p_0(  t+h-\tilde{p}^*_{1,1}(x))\nonumber\\
&\qquad-p_0(  t-\tilde{p}^*_{1,1}(x))\Big|p_1^*(x)\,dx\nonumber\\
&\quad=g(0)\lambda\int^{t}_{t-\tilde{p}^*_{1,1}(L)}\left|p_0(\tilde{t}+h)-p_0(\tilde{t})\right|\,d\tilde{t}.\label{eq: estimate of the failure state}
\end{align}
Moreover, in light of   \eqref{eq: estimate of the good state} we get
\begin{align}
&\int^L_0|p_1(x, t+h)-p_1(x, t)|\,dx  \nonumber\\
&\qquad= g(0)\lambda\int^{t}_{t-\tilde{p}^*_{1,1}(L)}\left|p_0(\tilde{t}+h)-p_0(\tilde{t})\right|\,d\tilde{t}\nonumber\\
&\qquad\leq2(g(0)\lambda)^2\tilde{p}^*_{1,1}(L) h\sup_{\tilde{t}\geq 0 } |p_0(\tilde{t})|. \label{eq: estimate of the failure state}
\end{align}
Thus, by  \eqref{eq: estimate of the good state} and~\eqref{eq: estimate of the failure state} we have
\begin{align*}
&\|T(t+h)P_0-T(t)P_0\|_{X}\\
&\leq (1+g(0)\lambda\tilde{p}^*_{1,1}(L))2g(0)\lambda h\sup_{t\geq 0 } |p_0(t)|\\&\leq C(h)\|P_0\|_X,
\end{align*}
where $C(h):=(1+g(0)\lambda\tilde{p}^*_{1,1}(L))2g(0)\lambda h$ and the last inequality follows from the fact that 
$\sup_{t\geq 0 } |p_0(t)|\leq  \sup_{t\geq 0 } \|T(t)P(0)\|_{X}$
 and $\|T(t)\|_{\mathcal{L}(X)}\leq 1$ for all $t\geq 0$.
Finally, since 
$
\lim_{h\rightarrow 0}C(h)=0,
$
we get
\[
\lim_{h\rightarrow 0}\|T(t+h)-T(t)\|_{X}=0,
\]
which completes the proof. \hfill$\Box$

According to \cite[p.\,331, Cor.\,3.3]{engel2000one} and \eqref{ss}, the following result holds immediately. 
\begin{corollary}\label{cor1}
For $P(0)\in X$, the time-dependent solution $P(x,t)= T(t)P(0)$ to \eqref{IVP_closed_loop} converges to its  steady-state solution $P^*(x)=(p^*_0, p^*_1(x))^T$ exponentially, that is,
\begin{align}
\Vert P(\cdot,t)-P^*(\cdot)\Vert_{X} \leq  M_{0}e^{-\varepsilon_{0}t},
\end{align}
for some constants $\varepsilon_{0}>0$ and $M_{0}\geq 1$.  
\end{corollary}

\section{Bilinear Controllability }
In this section, we present the proof of our main Theorem  \ref{main}.
Now consider the weighted closed-loop system  \eqref{1closed_loop}--\eqref{IC_F} for   $\alpha\geq p^*_1(0)$ and $i\in \mathbb{Z}^+$.

\noindent {\it Proof of Theorem \ref{main}.} 
Based on \eqref{1closed_loop}--\eqref{2closed_loop}, the closed-loop system is now weighted by $\alpha i$ for $t\in [t_{i-1}, t_i), i\in \mathbb{Z}^+$, and hence the decay rate of the system solution  to its steady-state
becomes $\alpha i\varepsilon_0$ for $t\in [t_{i-1}, t_i)$. Further note that $t_i-t_{i-1}=\frac{c_0}{i^2}$. Consequently, 
by Cor. \ref{cor1} we have
\begin{align}
&\|P(\cdot,c_0\sum^{i}_{k=1}\frac{1}{k^2} )-P^*(\cdot)\|_{X}\nonumber\\
&\quad\leq M_0e^{-c_0\sum^{i}_{k=1} \alpha  k \varepsilon_{0} \frac{1}{k^2} }
=M_0e^{-\alpha \varepsilon_{0}c_0 \sum^{i}_{k=1} \frac{1}{k}}. \label{EST_key}
\end{align}
 Since $t_f=c_0\sum^{\infty}_{k=1} \frac{1}{k^2}$ and $\lim_{i\to \infty}\sum^{i}_{k=1} \frac{1}{k}$ diverges, \eqref{EST_key} immediately implies
 $$ P(\cdot, t_f)=P^*(x),$$
which completes the proof.
\hfill$\Box$

\subsection{Boundedness of $\mu$ for $x\in [0, l]$ with $0<l<L$}
Finally,  if $p^*_{1}\in W^{1, \infty}(0, L)$, then one can   show that the  repair rate $\mu$ defined by  \eqref{mu_feedback} is bounded for $x\in[0, l]$ with $0<l<L$, 
when $t_f>2\|p^*_1\|_{L^1}$. Since $\frac{p_{1_x}}{p_1}$ converges to $\frac{p^*_{1_x}}{p^*_1}$ as $i\to \infty$, which is in $L^\infty(0, l)$ for $0<l<L$.
   It suffices to show that 
$\sup_{x\in [0, l]}|\alpha i g(x)p_{1}(x, t) |$ is finite  for $i\in \mathbb{Z}^+$ large enough. 
\begin{prop}
Let  $\alpha \geq \max\{p^*_1(0), \frac{1}{c_0}, \frac{1}{c_0\varepsilon_0}\}$.
For $t_f>2\|p^*_1\|_{L^1}$,
 the solution $p_1(x,t)$   to \eqref{1closed_loop}--\eqref{IC_F}  satisfies
\begin{align*}
&\lim_{i\to \infty}\sup_{x\in [0, l]}\left|\frac{\partial(\alpha i   g(x)p_1(x, t))}{\partial x}\right| <\infty.
\end{align*}
Further if $p^*_{1}\in W^{1\,\infty}(0, L)$, then the repair rate $\mu(x, t)$ defined by  \eqref{mu_feedback} is bounded for $x\in[0, l]$ with $0<l<L$.
\end{prop}

\noindent{\it Proof.}  
For simplicity, we denote  $(p_{0,i}, p_{1,i}(x,t))^T$ by $(p_{0}, p_{1}(x,t))^T$ in the rest of the proof.  For $t_f> 2\|p^*_1\|_{L^1}$, there exists $i\in\mathbb{Z}^+$ large enough such that $t_{i-1}>  2\|p^*_1\|_{L^1}$, where 
$\|p^*_1\|_{L^1}\geq \tilde{p}^*_{1,i}(L)$ for any $i\in \mathbb{Z}^+$ and $\alpha\geq1$. Thus for $t\geq t_{i-1}$ we have $t>   2\tilde{p}^*_{1,i}(L) $, 
and hence by \eqref{ex_sol_p1}, 
\begin{align*}
p_1(x, t)=\frac{1}{ g(x)} g(0)\lambda p_0(  t-\tilde{p}^*_{1,i}(x)).
\end{align*}
Therefore, 
\begin{align}
\frac{\partial (g(x)p_1(x,t))}{\partial x}=g(0)\lambda\frac{dp_0(  t-\tilde{p}^*_{1,i}(x))}{dt}(-\frac{p^*_1(x)}{\alpha i}).\label{EST_p1_1}
\end{align}
With the help of \eqref{ex_sol_p0},  we get
\begin{align}
\frac{dp_0(  t-\tilde{p}^*_{1,i}(x))}{dt} =&\alpha i g(0)\lambda\big(  p_0(  t-\tilde{p}^*_{1,i}(x)-\tilde{p}^*_{1,i}(L))\nonumber\\
&-p_0(  t-\tilde{p}^*_{1,i}(x))\big). \label{EST_p1_2}
\end{align}
Combining  \eqref{EST_p1_1} with \eqref{EST_p1_2} and \eqref{EST_key}  yields
\begin{align*}
& \sup_{x\in [0, l]} \left|\frac{\partial(\alpha i   g(x)p_1(x, t))}{\partial x}\right|\nonumber\\
&=\alpha i (g(0)\lambda)^2 \sup_{x\in [0, l]} \big|\big( p_0\big(t-\tilde{p}^*_{1,i}(x) -\tilde{p}^*_{1,i}(L)\big)\nonumber\\
&\quad-p_0\big(  t-\tilde{p}^*_{1,i}(x)\big)\big)p^*_1(x)\big| \nonumber\\
&\leq \alpha i (g(0)\lambda)^2  
  \sup_{x\in [0, l]}\big(\big|p_0( t-\tilde{p}^*_{1,i}(x)- \tilde{p}^*_{1,i}(L) )-p^*_0\big| \nonumber\\
&\quad+\big|p_0(  t-\tilde{p}^*_{1,i}(x)) -p^*_0\big|\big)\cdot \sup_{x\in [0, l]}| p^*_1(x) | \nonumber\\
&\leq \alpha^2 i (g(0)\lambda)^2 \nonumber\\
&\qquad\cdot\Big(  \sup_{\tau\in [t-\tilde{p}^*_{1,i}(l)- \tilde{p}^*_{1,i}(L) , t- \tilde{p}^*_{1,i}(L) ]}\big|p_0( \tau )-p^*_0\big| \nonumber\\
&\quad+\sup_{\tau\in [t-\tilde{p}^*_{1,i}(l), t]}\big|p_0( \tau) -p^*_0\big|\Big) \nonumber\\
&\leq 2 \alpha^2 i (g(0)\lambda)^2  \sup_{\tau\in [t-2 \tilde{p}^*_{1,i}(L) , t]}\big|p_0( \tau )-p^*_0\big| \nonumber\\
&\leq 2 \alpha^2 i (g(0)\lambda)^2 \|P(\cdot, t-2 \tilde{p}^*_{1,i}(L) )-P^*\|_{X}.%
\end{align*}

Recall that  $ \int^L_0p^*_1(x)\,dx\leq 1$ and  $\alpha\geq \frac{1}{c_0}$. We have 
 $$\tilde{p}^*_{1,i}(L)=\frac{1}{\alpha i}\int^L_0p^*_1(x)\,dx\leq \frac{c_0}{i}.$$
Then
 $ t-2\tilde{p}^*_{1,i}(L)\geq t_{i-1}-\frac{2c_0}{i}$ for $i\geq 2$.
In light of Corollary \ref{cor1} and \eqref{EST_key} we get  
\begin{align*}
 & \|P(\cdot, t-2 \tilde{p}^*_{1,i}(L) )-P^*\|_{X} \nonumber\\
 &\leq \|P(x, t_{i-1}-\frac{2c_0}{i})-P^*\|_{X} \nonumber\\
& \leq M_0e^{-\alpha  \varepsilon_{0} c_0\sum^{i-1}_{k=1} \frac{1}{k} }\cdot e^{\frac{2c_0}{i}\alpha \varepsilon_0(i-1) }. 
\end{align*}
Therefore, 
\begin{align}
&\left|\frac{\partial(\alpha i   g(x)p_1(x, t))}{\partial x}\right|\nonumber\\
 &\leq 2 \alpha^2 i  (g(0)\lambda)^2 M_0e^{-\alpha   \varepsilon_{0} c_0\sum^{i-1}_{k=1} \frac{1}{k} }\cdot e^{\frac{2c_0}{i}\alpha(i-1)  \varepsilon_0 },
 \label{EST_bd}
\end{align}
where 
\begin{align}
\lim_{i\to \infty}e^{\frac{2c_0}{i}\alpha \varepsilon_0 (i-1) }=e^{2c_0\alpha \varepsilon_0 }.\label{EST_bd2}
\end{align}
 It remains to analyze the property of 
$i  e^{-\alpha  \varepsilon_{0} c_0\sum^{i-1}_{k=1}   \frac{1}{k} }$ when $i$ is sufficiently large. 
Let $j=i-1$. Then 
\begin{align*}
&i  e^{-\alpha  \varepsilon_{0} c_0\sum^{i-1}_{k=1}   \frac{1}{k} } =(j+1)  e^{-\alpha  \varepsilon_{0} c_0\sum^{j}_{k=1}   \frac{1}{k} } \\
&= j e^{-\alpha  \varepsilon_{0} c_0\sum^{j}_{k=1}   \frac{1}{k} }+e^{-\alpha  \varepsilon_{0} c_0\sum^{j}_{k=1}   \frac{1}{k} }\\
&=e^{-(-\ln j+\alpha  \varepsilon_{0} c_0 \sum^{j}_{k=1}   \frac{1}{k} )}+ e^{-\alpha  \varepsilon_{0} c_0\sum^{j}_{k=1}   \frac{1}{k} }.
\end{align*}
Since 
$ \lim_{j\to\infty} (-\ln j+ \sum^{j}_{k=1}   \frac{1}{k} )=\gamma>0$
  is the Euler-Mascheroni constant \cite[Sec.1.5]{finch2003mathematical}, 
 we have for $\alpha\geq\frac{1}{ \varepsilon_{0} c_0}$,
 \begin{align}
 &\lim_{j\to \infty}(e^{-(-\ln j+\alpha  \varepsilon_{0} c_0 \sum^{j}_{k=1}   \frac{1}{k} )}+ e^{-\alpha  \varepsilon_{0} c_0\sum^{j}_{k=1}   \frac{1}{k} }) \leq e^{-\gamma}. \label{EST_bd3}
 \end{align}
Finally, combining  \eqref{EST_bd} with  \eqref{EST_bd2} and \eqref{EST_bd3} follows
$$\lim_{i\to \infty}\sup_{x\in [0, l]} \left| \frac{\partial(\alpha i   g(x)p_1(x, t))}{\partial x}\right|<\infty.$$
This completes the proof. \hfill$\Box$

\section{Conclusion }
Bilinear controllability of a simple reparable system via system repair rate is addressed in this work. A specific control law in feedback form is constructed. 
The construction essentially makes use of  the exponential convergence of the system solution to its steady-state. In fact, there are many other ways of choosing the control weight in  \eqref{mu_feedback}  as long as the series in \eqref{EST_key} diverges. Our approach is generic and can be applied to a broad family of reparable systems of similar attributes. 

\end{document}